\documentclass[11pt, a4]{amsart}
\usepackage[a4paper, left=28mm, right=28mm, top=28mm, bottom=34mm]{geometry}
\usepackage{amsmath} 					%数式環境
\usepackage{amsthm} 						%定理環境
\usepackage{amssymb} 					%記号。txfontsやpxfontsもあり
\usepackage{amscd} 						%簡単な図式。xypicが面倒なときに

\usepackage{epic, eepic}
\usepackage{graphicx} 		%図。上部下部にキャプション可能
\usepackage{here}					%表の位置の強制的な修正

\usepackage[all]{xy}							%可換図式ほか。可換図式のラベルに行列を用いるにはxymatrixなどの工夫が要る
\usepackage{longtable}

\usepackage{mathrsfs} 					%花文字
\usepackage{latexsym}						%証明終わりの記号のための何か
\usepackage{type1ec} % <-- New!
\usepackage[OT2,T1]{fontenc}
\usepackage[russian,english]{babel}
\usepackage{enumerate}					%番号つき箇条書き。descriptionでいい気はするが
\usepackage{verbatim}						%文字通りに書く環境。
\AtBeginDocument{\shorthandoff{"}}			
\newcommand{\Integer}{\mathbb{Z}}					%整数環
\def\det{\mathop{\mathrm{det}}\nolimits}			%行列式
\newcommand{\isomorphism}{\overset{\sim}{\to}}			%同型な射、ホモトピー同値
\def\Kernel{\mathop{\mathrm{Ker}}\nolimits}		%核
\def\sheafhom{\mathop{\mathscr{H}\kern -2pt om}\nolimits}		%層の射

\def\sheafend{\mathop{\mathscr{E}\kern -2pt nd}\nolimits}			%層の自己準同型
\def\sheafext{\mathop{\mathscr{E}\kern -2pt xt}\nolimits}			%層のExt
\def\Left{\mathop{\mathrm{L} \kern -2pt}\nolimits}				%左導来函手
\def\Right{\mathop{\mathrm{R} \kern -2pt}\nolimits}			%右導来函手
\newcommand{\Cohomology}[2]{H^{#1}\! \left( {#2} \right)}

\newcommand{\strshf}{\mathcal{O}}				%整数環・関数環
\newcommand{\projsp}{\mathbb{P}}				%射影空間。\Pは段落終わり記号
\def\Spec{\mathop{\mathrm{Spec}}\nolimits}			%スペクトル
\def\sheafspec{\mathop{\mathscr{S}\kern -2pt pec}\nolimits}
\def\sheafproj{\mathop{\mathscr{P}\kern -2pt roj}\nolimits}

\def\Pic{\mathop{\mathrm{Pic}}\nolimits}			%直線束の同値類
\def\Jac{\mathop{\mathrm{Jac}}\nolimits}		%Jacobians
\def\div{\mathop{\mathrm{div}}\nolimits}			%付随する因子
\def\Brauer{\mathop{\mathrm{Br}}\nolimits}		%Brauer群
\def\GenLin{\mathop{\mathrm{GL}}\nolimits}			%一般線型群
\def\SpecLin{\mathop{\mathrm{SL}}\nolimits}			%特殊線型群
\def\Matrix{\mathop{\mathrm{Mat}}\nolimits}		%行列環

\makeatletter
\def\iddots{\mathinner{\mkern1mu\raise\p@
    \hbox{.}\mkern2mu\raise4\p@\hbox{.}\mkern2mu
    \raise7\p@\vbox{\kern7\p@\hbox{.}}\mkern1mu}}
\def\adots{\mathinner{\mkern2mu\raise\p@\hbox{.} %% yhmath.styから
 \mkern2mu\raise4\p@\hbox{.}\mkern1mu
 \raise7\p@\vbox{\kern7\p@\hbox{.}}\mkern1mu}}
\makeatother

\theoremstyle{definition}
\newtheorem{thm}{Theorem}[section]			%定理
\newtheorem{prop}[thm]{Proposition}			%命題
\newtheorem{cor}[thm]{Corollary}					%系
\newtheorem{Alg}[thm]{Algorithm}			% Algorithm

\theoremstyle{definition}
\newtheorem*{prf}{Proof}								%証明
\newtheorem{rmk}[thm]{Remark}					%注意
 \makeatletter
    
    \@addtoreset{equation}{section}
  \makeatother

\newcommand{\relmiddle}[1]{\mathrel{}\middle#1\mathrel{}}	
\newcommand{\relmid}{\relmiddle{|}}	% ご所望の|

\newcommand{\Equref}[1]{(\ref{#1})}		%()付きの引用

\title
{Linear determinantal representations of smooth plane cubics 
over finite fields}
\date{\today}
\author{Yasuhiro Ishitsuka}
\address{Department of Mathematics, Faculty of Science, Kyoto University, Kyoto 606-8502, Japan}
\email{yasu-ishi@math.kyoto-u.ac.jp}


\begin{document}

\begin{abstract}
In this note, we study linear determinantal
representations of smooth plane cubics over finite fields.
We give an explicit formula of linear determinantal representations
corresponding to rational points. Using Schoof's formula, 
we count the number of projective equivalence classes of 
smooth plane cubics over a finite field admitting prescribed number of
equivalence classes of linear determinantal representations.
As an application, we determine isomorphism classes of 
smooth plane cubics over a finite field
with 0, 1 or 2 equivalence classes of linear determinantal representations.
\end{abstract}

\maketitle

\section{Introduction}
Let $k$ be a field, and
\begin{align*}
	F(X, Y, Z) = a_{000}X^3 &+ a_{001}X^2Y + a_{002}X^2Z +
	a_{011}XY^2 + a_{012}XYZ \\
	&+ a_{022}XZ^2 +
	a_{111}Y^3 + a_{112}Y^2Z + a_{122}YZ^2 + a_{222}Z^3
\end{align*}
a {ternary cubic form} with coefficients in $k$
defining a smooth plane cubic $C \subset \projsp^2$.
We say that the cubic \textit{$C$ admits
a linear determinantal representation over $k$} if
there are a nonzero constant $0 \neq \lambda \in k$ and
three square matrices $M_0, M_1, M_2 \in \Matrix_3(k)$
of size 3 satisfying
\(
	F(X, Y, Z) = \lambda \cdot \det (M),
\)
where we put $M := XM_0 + YM_1 + ZM_2$.
We say that two linear determinantal representations $M, M'$ of $C$ are
\textit{equivalent} if there are invertible matrices $A, B \in \GenLin_3(k)$
such that 
\(
	M' = AMB.
\)

Studying linear determinantal representations of smooth plane cubics is 
a classical topic in linear algebra and algebraic geometry 
(for example, see \cite{Vin89}, \cite{Dol12}).
Recently, they appear in the study of 
the derived category of smooth plane cubics
(\cite{Gal14}, \cite{BP15}), and have been studied from
arithmetic viewpoints (\cite{FN14}, \cite{II14}, \cite{Ish15}).

In this note, we investigate linear determinantal representations
of smooth plane cubics over finite fields.
Let $\mathbb{F}_q$ be a finite field with $q$ elements.
First, we prove the following bijection. 
Recall that any smooth plane cubic over $\mathbb{F}_q$ has 
a $\mathbb{F}_q$-rational point (\cite[Theorem 3]{Lan55}).
\begin{thm}[See Proposition \ref{Corr2} and 
Theorem \ref{Th: LDR2}]\label{main1}
	Let $C \subset \projsp^2$ be a smooth plane cubic over 
	$\mathbb{F}_q$. Fix an $\mathbb{F}_q$-rational
	point $P_0 \in C(\mathbb{F}_q)$.
	There is a natural bijection between the following two sets:
	\begin{itemize}
		\item the set of equivalence classes of 
		linear determinantal representations of $C$ over $\mathbb{F}_q$, and
		\item the set $C(\mathbb{F}_q) 
		\setminus \{P_0\}$ of $\mathbb{F}_q$-rational points on $C$
		different from $P_0$.
	\end{itemize}
\end{thm} 
We also calculate a representative of the equivalence class of
linear determinantal representations corresponding to
each $\mathbb{F}_q$-rational point $P \in C(\mathbb{F}_q) \setminus \{P_0\}$ 
(for a precise statement, see Theorem \ref{Th: LDR2}).
In fact, these results are valid for smooth plane cubics with rational points
over arbitrary fields.

%Our second main theorem treats the case when 
%$k$ is a finite field. Two cubics $C, C' \subset \projsp^2$ are
%said to be \textit{projectively equivalent} if there exists 
%an isomorphism $\rho \colon \projsp^2 \isomorphism \projsp^2$
%such that $\rho|_C$ gives an isomorphism $C \isomorphism C'$. 
Let $\mathrm{Cub}_q(n)$ be the number of projective equivalence classes 
of smooth plane cubics over $\mathbb{F}_q$ 
with exactly $n$ equivalence classes of linear determinantal 
representations. We compute $\mathrm{Cub}_q(n)$ for $0 \le n \le 2$.	

\begin{thm}[{See Corollary \ref{Th: FinLDR} and Section 
\ref{FinLDR}}]\label{main2}\mbox{}
	\begin{enumerate}
		\item For $2 \le q \le 4$, we have $\mathrm{Cub}_q(0)=1$;
		otherwise, $\mathrm{Cub}_q(0) = 0$.
		\item For $2 \le q \le 5$, we have $\mathrm{Cub}_q(1) = 1$;
		otherwise, $\mathrm{Cub}_q(1) = 0$.
		\item For $q = 2,3,5,7$, we have $\mathrm{Cub}_q(2) =2$.
		For $q=4$, we have $\mathrm{Cub}_q(2) = 4$.
		Otherwise, $\mathrm{Cub}_q(2) = 0$.
	\end{enumerate}
	
	\begin{table}[h]
		\begin{center}
			\caption{The number of projective equivalence classes of
			smooth plane cubics over finite fields admitting prescribed number of
			equivalence classes of linear determinantal representations.}
			\begin{tabular}{|c|c|c|c|c|c|c|}
				\hline
				 & $\mathbb{F}_2$ & $\mathbb{F}_3$ & 
				$\mathbb{F}_4$ & $\mathbb{F}_5$ & $\mathbb{F}_7$ & 
				$\mathbb{F}_q \; (q \ge 8)$
				\\ \hline
				$\mathrm{Cub}_q(0)$ & 1 & 1 & 1 & 0 & 0 & 0 
				\\ \hline
				$\mathrm{Cub}_q(1)$ & 1 & 1 & 1 & 1 & 0 & 0
				\\ \hline
				$\mathrm{Cub}_q(2)$ & 2 & 2 & 4 & 2 & 2 & 0
				\\ \hline
			\end{tabular}
			\label{Tb: FLDR0}
		\end{center}
	\end{table}
\end{thm}

For each equivalence class in this table,
we give examples of smooth plane cubics
and their linear determinantal representations.
In particular, we determine all projective equivalence classes of
smooth plane cubics over finite fields which admit at most
two equivalence classes of linear determinantal representations.
See Table \ref{Tb: 0LDR} to Table \ref{Tb: 2LDR-7}.

The outline of this paper is as follows.
In Section 2, we recall the notion of 
linear determinantal representations of smooth plane curves and
its relation to a class of line bundles.
In Section 3, we describe an algorithm to compute a representative of
linear determinantal representations corresponding to a line bundle.
Then we perform this algorithm to smooth plane cubics
with rational points, and obtain an explicit formula of
linear determinantal representations in Section 4.
In Section 5, we recall Schoof's formula counting
the number of projective equivalence classes of smooth plane cubics over
finite fields with prescribed number of rational points.
Then we apply it to count the number of projective equivalence classes of
smooth plane cubics over finite fields admitting prescribed number of
equivalence classes of linear determinantal representations.
Finally, in Section 6, we determine smooth plane cubics over finite fields 
admitting at most two equivalence classes of 
linear determinantal representations.

\section{Linear determinantal representations of smooth plane cubics 
with rational points}\label{RatLDR}
Let $k$ be a field, and $F(X, Y, Z) \in k[X, Y, Z]$ a homogeneous polynomial
with coefficients in $k$ of degree $d \ge 1$ 
defining a smooth plane curve $C \subset \projsp^2$. 
Its degree is $d$, and its genus is $g = (d-1)(d-2)/2$. 
We fix projective coordinates $X, Y, Z$ of $\projsp^2$.

	A \textit{linear determinantal representation} of $C$ over $k$ is
	a square matrix $M$ of size $d$ with entries in 
	$k$-linear forms in three variables $X, Y, Z$ which satisfies
		\(
			F(X, Y, Z) = \lambda \cdot \det(M)
		\)
	for some $\lambda \in k^\times$.
	Two linear determinantal representations $M, M'$ are 
	said to be \textit{equivalent} if there exist two invertible matrices
	$A, B \in \GenLin_d(k)$ with $M' = AMB$. 
	We denote by $\mathrm{LDR}(C)$ the set of equivalence classes of
	linear determinantal representations of $C$ over $k$.

The following theorem gives an interpretation of 
linear determinantal representations 
of $C$ in terms of non-effective line bundles on $C$.
It is well known at least when $k$ is 
an algebraically closed field of characteristic zero.

\begin{thm}[{see \cite[Proposition 3.1]{Bea00}, \cite[Proposition 2.2]{Ish15}}]\label{Corr}
	There is a natural bijection between the following two sets:
	\begin{itemize}
		\item the set $\mathrm{LDR}(C)$ of equivalence classes of 
		linear determinantal representations of $C$ over $k$, and
		\item the set of isomorphism classes of
		non-effective line bundles on $C$ of degree $g-1$.
	\end{itemize}
\end{thm}

\begin{prf}
	We briefly recall the proof because it is used to prove the correctness of
	the algorithm in Section 3. See also \cite{Bea00}, \cite{Ish14}, \cite{Ish15}
	for details.
	
	We take a non-effective line bundle $\mathcal{L}$ of degree $g-1$ on $C$.
	Let $\iota \colon C \hookrightarrow \projsp^2$ be the given embedding.
	We denote the homogeneous coordinate ring of $\projsp^2$ by
	\begin{align*}
		R &:= \Gamma_*(\projsp^2, \strshf_{\projsp^2})\\
		&= \bigoplus_{n \in \Integer} \Cohomology{0}{\projsp^2,
		\strshf_{\projsp^2}(n)}.
	\end{align*}
	The graded $R$-module $N = \Gamma_*(\projsp^2, \iota_*\mathcal{L})
	\cong \Gamma_*(C, \mathcal{L})$ has 
	a minimal free resolution of the form
	\begin{equation}\label{Eq: LDR}
		\xymatrix{
			0 \ar[r] & R(-2) \otimes_k W_1 \ar[r]^{\widetilde{M}} & 
			R(-1) \otimes_k W_0 \ar[r] & N \ar[r] & 0
		},
	\end{equation}
	where $W_0, W_1$ are $d$-dimensional $k$-vector spaces \cite[Proposition 3.1]{Bea00}.
	The homomorphism $\widetilde{M}$ can be 
	expressed by a square matrix $M$ of size $d$
	with coefficients in $k$-linear forms in three variables $X, Y, Z$.
	We can check $M$ gives a linear determinantal representation of $C$,
	and its equivalence class depends only on 
	the isomorphism class of the line bundle $\mathcal{L}$.
	
	Conversely, we take a linear determinantal representation $M$ of $C$.
	This matrix gives an injective homomorphism 
	\[
		\widetilde{M} \colon R(-2)^{\oplus d} \to R(-1)^{\oplus d}. 
	\]
	We denote by $N$ the cokernel of $\widetilde{M}$.
	We can show that the coherent sheaf associated to $N$ is 
	written as $\iota_* \mathcal{L}$
	for a non-effective line bundle $\mathcal{L}$ 
	of degree $g-1$ on $C$. 
	The isomorphism class of $\mathcal{L}$ depends only on
	the equivalence class of $M$. 	
	By construction, these two maps are inverses to each other.
\qed\end{prf}

Assume that $d=3$, i.e., $C$ is a smooth plane cubic over $k$.
We shall study the relation between the Picard group $\Pic^0(C)$ and
the group $\Jac(C)(k)$ of $k$-rational points on 
the Jacobian variety $\Jac(C)$ of $C$.
In general, there can be a difference which is measured by
the relative Brauer group (for example, see \cite[Theorem 2.1]{CK12}, 
\cite[Example 6.9]{Ish15}).
However, when $C$ has a $k$-rational point, the difference vanishes.

\begin{prop}\label{Corr2}
	Let $C$ be a smooth plane cubic over $k$ 
	with a $k$-rational point $P_0 \in C(k)$.
	There is a natural bijection between the following two sets:
	\begin{itemize}
		\item the set $\mathrm{LDR}(C)$ of equivalence classes of 
		linear determinantal representations of $C$ over $k$, and
		\item the set $C(k) \setminus \{P_0\}$ 
		of $k$-rational points on $C$ different from $P_0$.
	\end{itemize}
\end{prop}

\begin{prf}
	There is an exact sequence
	\[
		\xymatrix{
			0 \ar[r] & \Pic(C) \ar[r] & \Pic_{C/k}(k) \ar[r] &
			\Brauer(k) \ar[r]^s & \Brauer(C),
		}
	\]
	where $s$ is the pullback morphism associated to the structure morphism
	$C \to \Spec(k)$ (\cite[Theorem 2.1]{CK12}). Since $C$ has a $k$-rational point,
	the homomorphism $s$ is injective. Hence we have two isomorphisms
	\[
		\Pic(C) \isomorphism \Pic_{C/k}(k), 
		\quad \Pic^0(C) \isomorphism \Jac(C)(k).
	\]
	Then the morphism
	\begin{align*}
		\iota_{P_0} \colon C &\to \Jac(C)\\
		P &\mapsto P - P_0
	\end{align*}
	gives an isomorphism. 
	The only effective line bundle on $C$ of degree 0 is
	the trivial bundle $\strshf_C = \iota_{P_0}(P_0)$. 
	Thus, by Theorem \ref{Corr} and the bijection $\iota_{P_0}$, 
	we have the desired bijection.
\qed\end{prf}

\section{An algorithm to obtain linear 
determinantal representations of smooth plane curves}

Let us make the bijection in Theorem \ref{Corr} explicit.
In this section, we shall give an algorithm 
to obtain a linear determinantal representation
of a smooth plane curve $C$ of degree $d$ and
genus $g = (d-1)(d-2)/2$ over an arbitrary field $k$.
\begin{Alg}\label{Alg:LDR}\mbox{} \\[-13pt]
\begin{description}
	\item[Input] a defining equation $F(X,Y,Z)$ of $C \subset \projsp^2$ 
	with respect to fixed projective coordinates $X, Y, Z$, and
	a $k$-rational non-effective divisor $D$ of degree $g-1$.
	\item[Output] a linear determinantal representation of $C$ over $k$ 
	corresponding to $D$. 
	\begin{description}
		\item[Step 1 (Global Section)] Compute a $k$-basis $\{v_0, v_1, v_2\}$
		of the 3-dimensional $k$-vector space
		$\Cohomology{0}{C, \strshf_C(D)(1)}$.
		\item[Step 2 (First Syzygy)] Compute a $k$-basis $\{e_0, e_1, e_2\}$
		of the 3-dimensional $k$-vector space
		\[
			\Kernel \left(
				\Cohomology{0}{C, \strshf_C(1)} \otimes_k
				\Cohomology{0}{C, \strshf_C(D)(1)} \to
				\Cohomology{0}{C, \strshf_C(D)(2)}
			\right).
		\]
		\item[Step 3 (Output Matrix)] 
		Write the $k$-basis $\{e_0, e_1, e_2\}$ as
		\[
			e_i = \sum_{j} l_{i, j}(X, Y, Z) \otimes v_j,
		\]
		where $l_{i,j}(X, Y, Z) \in \Cohomology{0}{C, \strshf_C(1)}$ 
		are $k$-linear forms.
		Output the matrix \[M = (l_{i,j}(X, Y, Z)).\]
	\end{description}
\end{description}
\end{Alg}

\begin{prf}[Proof of the correctness of Algorithm \ref{Alg:LDR}]
	Recall the short exact sequence \Equref{Eq: LDR}
	\begin{equation*}\label{Eq: LDR2}
		\xymatrix{
			0 \ar[r] & R(-2) \otimes_k W_1 \ar[r]^{\widetilde{M}} & 
			R(-1) \otimes_k W_0 \ar[r] & N \ar[r] & 0,
		}
	\end{equation*}
	where $W_0, W_1$ are 3-dimensional $k$-vector spaces.
	Since $R_0 = k$ and $N = \Gamma_*(C, \strshf_C(D))$ is the graded $R$-module
	corresponding to $\strshf_C(D)$, the degree 1 part of this sequence gives
	\[
		W_0 = N_1 = \Gamma(C, \strshf_C(D)(1)).
	\]
	The degree 2 part gives a short exact sequence
	\begin{equation*}\label{Eq: LDR3}
		\xymatrix{
			0 \ar[r] & W_1 \ar[r]^(.4){\widetilde{M}} & 
			R_1 \otimes_k W_0 \ar[r] & N_2 \ar[r] & 0.
		}
	\end{equation*}
	Thus we have
	\[
		W_1 = \Kernel \left(
				\Cohomology{0}{C, \strshf_C(1)} \otimes_k
				\Cohomology{0}{C, \strshf_C(D)(1)} \to
				\Cohomology{0}{C, \strshf_C(D)(2)}
			\right).
	\]
	The morphism $\widetilde{M}$ is the canonical embedding
	$W_1 \to R_1 \otimes_k W_0$.
	Hence it is represented by
	the matrix $M = (l_{i,j}(X, Y, Z))$.
\qed\end{prf}

\section{An explicit formula on linear determinantal representations
of smooth plane cubics with rational points}
We apply Algorithm \ref{Alg:LDR} 
to a smooth plane cubic (i.e., $d=3$) with a $k$-rational point.
Note that, by changing projective coordinates,
we may assume that the smooth plane cubic $C$ over $k$ has
a $k$-rational point $P_0 = [1:0:0]$, and 
the tangent line of $C$ at $P_0$ is $(Z=0)$.

\begin{thm}\label{Th: LDR2}
	Let $C \subset \projsp^2$ be a smooth plane cubic 
	over an arbitrary field $k$
	with a $k$-rational point $P_0 = [1:0:0]$.
	Assume that the tangent line of $C$ at $P_0$ is the line $l=(Z=0)$.
	We have the following formula for the equivalence class of
	linear determinantal representations of $C$ over $k$
	corresponding to a point $P = [s:t:u] \in C(k) \setminus \{P_0\}$ 
	via Proposition \ref{Corr2}.
	\begin{description}
		\item[Case 1] 
		If $u \neq 0$, the equivalence class of 
		linear determinantal representations of $C$
		corresponding to $P$ is given by
		\begin{align}\label{DRformula1}
			M_P = \begin{pmatrix}
				0 & Z & -Y \\
				uY - tZ & 0 & -u^2 X - (Q(t, u) + su) Z \\
				u X - s Z & L_1(X, Y, Z) & L_2(X, Y, Z)
			\end{pmatrix},
		\end{align}
		where we denote
		\begin{align*}
			L_1(X, Y, Z) &:= u^2 a_{011} X + u^2 a_{111} Y 
			+u (a_{111}t + a_{112}u) Z,\\
			L_2(X, Y, Z) &:=u(a_{011}t + a_{012}u) X
			+ (a_{111}t^2 + a_{112}tu + a_{122}u^2)Z, \\
			Q(Y, Z) &:= a_{011}Y^2 + a_{012}YZ + a_{022}Z^2.
		\end{align*}
		\item[Case 2] 
		If $u=0$, the equivalence class of 
		linear determinantal representations of $C$
		corresponding to $P$ is given by
		\begin{align}\label{DRformula2}
			M_P = \begin{pmatrix}
				0 & Z & -Y \\
				Z & a_{011}Y &  X + a_{012}Y + a_{022} Z \\
				a_{011}X + a_{111}Y & \widetilde{L}_1(X, Y, Z) 
				& \widetilde{L}_2(X, Y, Z)
			\end{pmatrix},
		\end{align}
		where we denote
		\begin{align*}
			\widetilde{L}_1(X, Y, Z) &:= a_{111}X 
			+(a_{012}a_{111} - a_{011}a_{112})Y,\\
			\widetilde{L}_2(X, Y, Z) &:= (a_{022}a_{111} - a_{011}a_{122})Y
			- a_{011}a_{222}Z.
		\end{align*}
	\end{description} 
\end{thm}

We shall prove Theorem \ref{Th: LDR2} by performing Algorithm \ref{Alg:LDR} as follows.

\subsection{Preparation}
By the condition of Theorem \ref{Th: LDR2},
we may assume that $a_{000} = a_{001} = 0$ and $a_{002} = 1$.
Thus we can take a defining equation of 
the given cubic $C \subset \projsp^2$ as
\[
	ZX^2 + Q(Y, Z)X + C(Y, Z) = 0,
\]
where $Q(Y,Z)$ is a binary quadratic form defined in the statement of 
Case 1 of Theorem \ref{Th: LDR2}, and we denote
\begin{align*}
	C(Y, Z) &:= a_{111} Y^3 + a_{112} Y^2Z +
	a_{122} YZ^2 + a_{222} Z^3.
\end{align*}
The divisor $l \cap C$ on $C$ can be written as $2P_0 + R$, where
\[
	R = [a_{111} : -a_{011} : 0 ].
\] 
Note that $R$ may or may not be equal to $P_0 = [1:0:0]$.

Take a point $P = [s: t: u] \in C(k) \setminus \{P_0\}$.
The line $m = \overline{PP_0}$ is defined by 
\[
	m(Y, Z) := uY - tZ.
\]
The divisor $m \cap C$ on $C$ is $P + P_0 + S$, where
\begin{align*}
	S &= [Q(t,u) +su : -tu: -u^2] \in C(k).
\end{align*}
Since $P - P_0 = \div(m) - 2 P_0 - S$,
the $k$-vector space $W_0 = \Gamma(C, 
\strshf_C(P - P_0)(1))$ is isomorphic to the $k$-vector space 
\[
	V=\left\{ {q(X, Y, Z)} \relmid 
	\begin{array}{l}
		q(X, Y, Z) \in \Gamma(X, \strshf_C(2)), \\
		\div q(X, Y, Z) - 2P_0 - S \ge 0
	\end{array}
	\right\}
\]
via the isomorphism
\(
	W_0 \to V ; f \mapsto fm.
\)
Consider a $k$-basis $\{ X^2, XY, Y^2, XZ, YZ, Z^2\}$ 
of $\Gamma(C, \strshf_C(2))$.
The first two elements $X^2, XY$ have order 0, 1 at $P_0 \in C(k)$,
and the other elements $XZ, Y^2, YZ, Z^2$ have order not less than 2
at $P_0 \in C(k)$. Hence for a quadratic form $q \in V$, 
we can write the quadratic form $q(X, Y, Z)$ as
\[
	q(X, Y, Z) = b_{02}XZ + b_{11}Y^2 + b_{12}YZ
	+ b_{22}Z^2
\]
for some constants $b_{02}, b_{11}, b_{12}, b_{22} \in k$ and
$q$ vanishes at $S$.
We divide the proof of Theorem \ref{Th: LDR2} into
two cases described in the statement: $u \neq 0$ and $u = 0$.

\subsection{Proof of Case 1: when $u \neq 0$}
In this case, we see that $l \neq m$ and $S \neq P_0$.
When $\div q(X, Y, Z) - 2P_0 -S \ge 0$, we have
\[
	u^2( -b_{02}(Q(t,u)+su) + b_{11}t^2 +
	b_{12}tu + b_{22}u^2) = 0.
\]
We can take a $k$-basis of $W_0$ as
\begin{align*}
	v_0 &:= (u^2XZ + (Q(t, u) + su) Z^2)/m, \\
	v_1 &:= Y, \\
	v_2 &:= Z.
\end{align*}
Next we compute a $k$-basis of the first syzygy module 
\[
	W_1 =
	\Kernel \left( \Gamma(C, \strshf_C(1)) \otimes_k W_0 \to
	\Gamma(C, \strshf_C(2)) \right).
\]
We find
\begin{align*}
	e_0 &= Z \otimes v_1 - Y \otimes v_2,\\
	e_1 &= (uY - tZ) \otimes v_0 - (u^2X + (Q(t, u) + su) Z) \otimes v_2, \\
	e_2 &= (u X - s Z) \otimes v_0 + 
	L_1(X, Y, Z) \otimes v_1 + L_2(X, Y, Z) \otimes v_2 
\end{align*}
form a $k$-basis of the first syzygy module $W_1$,
where $L_1(X, Y, Z), L_2(X, Y, Z)$ are linear forms defined 
in the statement of Theorem \ref{Th: LDR2}.
The corresponding determinantal representation is
\begin{align}
	M_P = \begin{pmatrix}
		0 & Z & -Y \\
		uY - tZ & 0 & - u^2 X - (Q(t, u) + su) Z \\
		u X - s Z & L_1(X, Y, Z) & L_2(X, Y, Z)
	\end{pmatrix}.
\end{align}
We may check that $\det(M_P) = -u^3 f$.
This proves Case 1 of Theorem \ref{Th: LDR2}.

\subsection{Proof of Case 2: when $u = 0$}
In this case, $S = P_0 = [1:0:0]$ and $l = m$.
We can take a $k$-basis of $W_0$ as
\begin{align*}
	v_0 &:= -(XZ + a_{011}Y^2 + a_{012}YZ 
	+ a_{022}Z^2)/Z,\\
	v_1 &:= Y, \\
	v_2 &:= Z.
\end{align*}

Next we compute a $k$-basis of the first syzygy module $W_1$.
We find
\begin{align*}
	e_0 &= Z \otimes v_1 - Y \otimes v_2,\\
	e_1 &= Z \otimes v_0 + a_{011}Y \otimes v_1 
	+ (X + a_{012}Y + a_{022}Z) \otimes v_2, \\
	e_2 &= (a_{011}X + a_{111}Y) \otimes v_0 + 
	\widetilde{L}_1(X, Y, Z) \otimes v_1 
	+ \widetilde{L}_2(X, Y, Z) \otimes v_2 
\end{align*}
form a $k$-basis of $W_1$,
where $\widetilde{L}_1(X, Y, Z), \widetilde{L}_2(X, Y, Z)$ are 
$k$-linear forms defined in the statement of Theorem \ref{Th: LDR2}.
The corresponding linear determinantal representation is
\begin{align}
	M_P = \begin{pmatrix}
		0 & Z & -Y \\
		Z & a_{011}Y & X + a_{012}Y + a_{022} Z \\
		a_{011}X + a_{111}Y & \widetilde{L}_1(X, Y, Z) 
		& \widetilde{L}_2(X, Y, Z)
	\end{pmatrix}.
\end{align}
We may check that $\det(M_P) = a_{011}f$.
This proves Case 2 of Theorem \ref{Th: LDR2}.\qed

\begin{rmk}\label{Rm: Other1}
	Let $k$ be a field of characteristic not equal to 2 nor 3, and 
	\begin{align}\label{Eq: Wform}
		E \colon ( Y^2Z - X^3 - aXZ^2 - b Z^3 = 0)
		\subset \projsp^2
	\end{align}
	an elliptic curve over $k$ with origin $P_0 = [0:1:0]$
	defined by a Weierstrass equation.
	Let $P=[\lambda: \mu: 1] \in E(k)$ be a $k$-rational point on 
	an affine part of $E$.
	Galinat gave in \cite[Lemma 2.9]{Gal14}
	a representative of linear determinantal representations of $E$ over $k$
	corresponding to the divisor $P-P_0$ of degree 0 as
	\[
		M'_P := \begin{pmatrix}
			X - \lambda Z & 0 & -Y-\mu Z \\
			\mu Z - Y & X + \lambda Z & (a + \lambda^2) Z\\
			0 & Z & -X
		\end{pmatrix}.
	\]
	Theorem \ref{Th: LDR2} gives an essentially same representative
	of linear determinantal representation in this case;
	actually, we can transform $M_P$ into $M'_P$ 
	by changing coordinates and elementary transformation.
	When $k$ is algebraically closed, Vinnikov \cite{Vin89}
	gave other representatives.
\end{rmk}

\begin{rmk}\label{Rm: Other2}
	Let $k$ be a field of characteristic not equal to 2 nor 3, and 
	\[
		C \colon ( X^3 + Y^3 + Z^3 + \lambda XYZ = 0)
		\subset \projsp^2
	\]
	a smooth plane cubic over $k$ defined by Hesse's normal form.
	Let 
	\(
		P = [a_0 : a_1: a_2] \in C(k)
	\) 
	be a $k$-rational point with $a_0a_1a_2 \neq 0$. 
	In \cite[Theorem A]{BP15}, 
	Buchweitz and Pavlov showed that the Moore matrix
	\[
		M''_P := \begin{pmatrix}
			a_0 X & a_1 Z & a_2 Y \\
			a_1 Y & a_2 X & a_0 Z\\
			a_2 Z & a_0 Y & a_1 X
		\end{pmatrix}
	\]
	gives a linear determinantal representation of $C$ over $k$
	corresponding to the divisor $3P - H$ of degree 0,
	where $H$ is a hyperplane section of $C$.
	Note that, when $k$ is not algebraically closed, 
	there can be a linear determinantal representation of $C$
	over $k$ which is not equivalent to any Moore matrices.
	Also the Moore matrices of two distinct $k$-rational points $P, P' \in C(k)$ 
	can give equivalent linear determinantal representations of $C$ over $k$.
	These are explained by the fact that the homomorphism
	\begin{align*}
		C = \Pic^1(C) &\to \Pic^3(C) \cong \Pic^0(C)\\
		P \hspace{14pt} &\mapsto 
		\hspace{11pt} 3P \hspace{11pt} \mapsto 3P - H
	\end{align*}
	is not an isomorphism in general.
\end{rmk}

\begin{rmk}
	To compute the Cassels--Tate pairing on the 3-Selmer groups 
	of elliptic curve,
	Fisher and Newton \cite{FN14} considered linear determinantal
	representations when $k$ is a number field, and
	$C$ is locally soluble but \textit{has no $k$-rational point}.
\end{rmk}

\section{A counting on smooth plane cubics over finite fields}
Let $p$ be a prime number, and $m \ge 1$ a positive integer.
Let $\mathbb{F}_q$ be a finite field with $q = p^m$ elements.
We recall Schoof's formula on the number of
the projective equivalence classes of 
smooth plane cubics over $\mathbb{F}_q$ 
with prescribed number of $\mathbb{F}_q$-rational points.
Here, two smooth plane cubics $C, C' \subset \projsp^2$ 
over $\mathbb{F}_q$ 
are said to be \textit{projectively equivalent} if there exists an isomorphism
$\projsp^2 \isomorphism \projsp^2$ over $\mathbb{F}_q$ that induces
an isomorphism $C \isomorphism C'$.

\begin{thm}[{\cite[Theorem 5.2]{Sch87}}]\label{Th: LDRnum}
	For an integer $n \in \Integer$, 
	the number of projective equivalence classes of
	smooth plane cubics $C$ over $\mathbb{F}_q$ 
	with $\#C(\mathbb{F}_q)=n$ is
	\begin{equation}\label{Eq: SchCnt}
		\#E_q(n) + \# E_{q,3}(n) + 3 \# E_{q,3,3}(n) 
		- \varepsilon_q(q + 1-n).
	\end{equation}
\end{thm}
Here, we use the following notation which is 
slightly different from \cite{Sch87}.
For reader's convenience, we recall the definition and formulas
for the terms appearing in \eqref{Eq: SchCnt}.
\begin{itemize}
	\item For an integer $a \in \Integer$ and a prime number $p$,
	$\left(a / p\right)$ denotes the Jacobi symbol.
	\item For a negative integer $\Delta \in \Integer_{<0}$ with 
	$\Delta \equiv 0, 1 \pmod 4$,
	\textit{Kronecker's class number} $H(\Delta)$ is defined to be
	the number of $\SpecLin_2(\Integer)$-orbits of 
	positive definite integral binary quadratic forms 
	\[
		\left\{ f(U, V) = aU^2 + bUV + cV^2 \in \Integer[U, V]
		\relmid a>0, b^2-4ac = \Delta
		\right\}
	\]
	with discriminant $\Delta$.
	Here $\gamma = \begin{pmatrix} 
	p & q\\
	r & s
	\end{pmatrix} \in \SpecLin_2(\Integer)$ acts on $f(U, V)$ as
	\[
		(\gamma \circ f)(U, V) = a(pU+rV)^2 + b(pU+rV)(qU+sV)
		+c(qU+sV)^2.
	\]
	\item Let $E_q(n)$ denote the set of isomorphism classes of
	elliptic curves over $\mathbb{F}_q$ with $\#E(\mathbb{F}_q) = n$.
	(In \cite{Sch87}, Schoof used $N(q+1-n)$ instead of $\#E_q(n)$.)
	From \cite[Theorem 4.6]{Sch87}, we have the following formula.
	\begin{itemize}
		\item If $t^2 > 4q$, we have $\#E_q(q + 1 - t) = 0$.
		\item If $t^2 \le 4q$ and $p \nmid t$, we have
			$\#E_q(q+1-t) = H(t^2-4q)$.
		\item If $t^2 \le 4q$, $t \equiv 0 \pmod p$ 
		and $m \equiv 1 \pmod 2$, the case is divided into three cases:
		\begin{itemize}
			\item If $t=0$, we have $\#E_q(q+1-t) = H(-4p)$.
			\item If $(t^2, p) = (2q, 2)$ or $(3q, 3)$, we have
			$\#E_q(q+1-t) = 1$.
			\item Otherwise, we have $\#E_q(q+1-t) = 0$.
		\end{itemize}		
		\item If $t^2 \le 4q$, $t \equiv 0 \pmod p$ 
		and $m \equiv 0 \pmod 2$, the case is divided into four cases:
		\begin{itemize}
			\item If $t=0$, we have $\#E_q(q+1-t) = 1-(-4/p)$.
			\item If $t^2=q$, we have $\#E_q(q+1-t) = 1-(-3/p)$.
			\item If $t^2=4q$, we have 
			\[
				\#E_q(q+1-t) = \frac{1}{12}(p + 6 - 4(-3/p) - 3(-4/p)).
			\]
			\item Otherwise, we have $\#E_q(q+1-t) = 0$.
		\end{itemize}
	\end{itemize}
	\item Let $E_{q,3}(n)$ denote the set of isomorphism classes of
	elliptic curves $E \in E_q(n)$ with non-trivial 3-torsion points.
	(In \cite{Sch87}, Schoof used $N_3(q+1-n)$ instead of $E_{q,3}(n)$.)
	It is easily described as
	\[
		E_{q, 3}(n) = \begin{cases}
			E_{q}(n) & (3 \mid n) \\
			\emptyset & (3 \nmid n).
		\end{cases}
	\]
	\item Let $E_{q,3,3}(n)$ denote the set of isomorphism classes of
	elliptic curves $E \in E_q(n)$ with 
	\[
		E(\mathbb{F}_q)[3] \cong (\Integer / 3 \Integer) ^2.
	\]
	(In \cite{Sch87}, Schoof used $N_{3\times 3}(q+1-n)$ 
	instead of $E_{q,3,3}(n)$.)
	From \cite[Theorem 4.9]{Sch87}, we have the following formula.
	\begin{itemize}
		\item We assume that the following four conditions are satisfied:
		$q \equiv 1 \pmod 3$, $t^2 \le 4q$, $p \nmid t$ and 
		$t \equiv q+1 \pmod 9$.
		Then we have
		\[
			\#E_{q,3,3}(q+1-t)= H\left( \dfrac{1}{9} (t^2-4q) \right).
		\]
		\item We assume that the following three conditions are satisfied:
		$2 \mid m$,
		$p \neq 3$ and
		$t = 2 \cdot (p/3)^{m/2} \cdot p^{m/2}$.
		Then we have
		\(
			\#E_{q,3,3}(q+1-t)= \#E_q(q+1-t).
		\)
		\item Otherwise, we have
		\(
			\#E_{q,3,3}(q+1-t)=0.
		\)
	\end{itemize}
	\item We set $t_0 \in \Integer \cup \{\infty\}$ as follows.
	\begin{enumerate}
		\item If $q \not\equiv 1 \pmod 3$, then we set
		$t_0 := \infty$. Note that, in this case, we always have $t \neq t_0$.
		\item If $p \not\equiv 1 \pmod 3$ but $q \equiv 1 \pmod 3$,
		we set
		\(
			t_0 :=2 \cdot \left( p/3 \right)^{m/2} \cdot p^{m/2}
		\)
		(note that, in this case, $m$ is even).
		\item If $p \equiv 1 \pmod 3$, 
		$t_0$ is the unique integer satisfying
		$t \equiv q + 1 \pmod 9$,
		$p \nmid t$ and
		$t^2 + 3x^2 = 4q$ for some integer $x \in \Integer$.
	\end{enumerate}
	\item We set $t_1 \in \Integer \cup \{\infty\}$ as follows.
	\begin{enumerate}
		\item If $q \not\equiv 1 $ nor $4 \pmod {12}$, then we set
		$t_1 := \infty$. Note that, in this case, we always have $t \neq t_1$.
		\item If $p \not\equiv 1 \pmod 4$ but $q \equiv 1 $ or $4 \pmod {12}$,
		we set
		\(
			t_1 := 2 \cdot \left( p/3 \right)^{m/2} \cdot p^{m/2}
		\)
		(note that, in this case, $m$ is even).
		\item If $p \equiv 1 \pmod 4$ and $q \equiv 1 $ or $4 \pmod {12}$, 
		$t_1$ is the integer satisfying $t \equiv q + 1 \pmod 9$, $p \nmid t$ and
		$t^2 + 4x^2 = 4q$ for some integer $x \in \Integer$.
	\end{enumerate}
	\item We define a function $\varepsilon_q(t)$ as follows:
	\[
		\varepsilon_q(t) := \begin{cases}
			2 & (t \in \{t_0, t_1\}, \mbox{ but } t_0 \neq t_1) \\
			3 & (t = t_0 = t_1 \mbox{ and } p = 2) \\
			4 & (t = t_0 = t_1 \mbox{ and } p \neq 2) \\
			0 & (\mbox{otherwise}).
		\end{cases}
	\]
\end{itemize}

By Proposition \ref{Corr2} and Theorem \ref{Th: LDRnum},
we have the following corollary.

\begin{cor}\label{Th: FinLDR}
	With the above notation, the number $\mathrm{Cub}_q(n)$ of 
	projective equivalence classes of
	smooth plane cubics $C$ over $\mathbb{F}_q$ 
	with $\# \mathrm{LDR}(C) = n$ is
	\begin{equation}\label{Eq:LDRFmla}
		\mathrm{Cub}_q(n) = \#E_q(n + 1) + 
		\# E_{q,3}(n + 1) + 3 \# E_{q,3,3}(n + 1) 
		- \varepsilon_q(q -n).
	\end{equation}
\end{cor}
\begin{prf}
	By Proposition \ref{Corr2}, we have $\#\mathrm{LDR}(C)=\#C(\mathbb{F}_q) -1$
	for a smooth plane cubic $C$ over $\mathbb{F}_q$.
	Using this and Theorem \ref{Th: LDRnum},
	we have the desired result. 
\qed\end{prf}

\begin{rmk}
The following table summarizes the values of $\mathrm{Cub}_q(n)$ for small $n$.
\begin{table}[h]
	\begin{center}
		\caption{The number of projective equivalence classes of
		smooth plane cubics over finite fields admitting prescribed number of
		equivalence classes of linear determinantal representations.}
		\begin{tabular}{|c|c|c|c|c|c|c|}
			\hline
			 & $\mathbb{F}_2$ & $\mathbb{F}_3$ & 
			$\mathbb{F}_4$ & $\mathbb{F}_5$ & $\mathbb{F}_7$ & 
			$\mathbb{F}_q \; (q \ge 8)$
			\\ \hline
			$\mathrm{Cub}_q(0)$ & 1 & 1 & 1 & 0 & 0 & 0 
			\\ \hline
			$\mathrm{Cub}_q(1)$ & 1 & 1 & 1 & 1 & 0 & 0
			\\ \hline
			$\mathrm{Cub}_q(2)$ & 2 & 2 & 4 & 2 & 2 & 0
			\\ \hline
		\end{tabular}
		\label{Tb: FLDR1}
	\end{center}
\end{table}\\
To check this, Table \ref{Tb: FLDR1-help} is helpful.
\begin{table}[h]
	\begin{center}
		\caption{The numbers appearing in the formula \eqref{Eq:LDRFmla}
		for $0 \le n \le 2$ and $2 \le q \le 7$.}
		\begin{tabular}{|c|c|c|c|c|c|c|}
			\hline
			& $\# E_q(1)$ & $\# E_q(2)$ & $\# E_q(3)$ &
			$\# E_{q,3}(1)$ & $\# E_{q,3}(2)$ & $\# E_{q,3}(3)$  
			\\ \hline
			$\mathbb{F}_2$ & 1 & 1 & 1 & 0 & 0 & 1
			\\ \hline
			$\mathbb{F}_3$ & 1 & 1 & 1 & 0 & 0 & 1
			\\ \hline
			$\mathbb{F}_4$ & 1 & 1 & 2 & 0 & 0 & 2
			\\ \hline
			$\mathbb{F}_5$ & 0 & 1 & 1 & 0 & 0 & 1
			\\ \hline
			$\mathbb{F}_7$ & 0 & 0 & 1 & 0 & 0 & 1
			\\ \hline 
		\end{tabular}
		\end{center}
%\end{table}
%\vspace{-10pt}
%\begin{table}[h]
	\begin{center}
		\begin{tabular}{|c|c|c|c|c|c|c|c|c|c|}
			\hline
			& $\# E_{q,3,3}(1)$ & $\# E_{q,3,3}(2)$ & $\# E_{q,3,3}(3)$ &
			$t_0$ & $t_1$ & $\varepsilon_{q}(q)$ 
			& $\varepsilon_{q}(q-1)$ & $\varepsilon_{q}(q-2)$  
			\\ \hline
			$\mathbb{F}_2$ & 0 & 0 & 0 & $\infty$ & $\infty$ & 0 & 0 & 0
			\\ \hline
			$\mathbb{F}_3$ & 0 & 0 & 0 & $\infty$ & $\infty$ & 0 & 0 & 0
			\\ \hline
			$\mathbb{F}_4$ & 0 & 0 & 0 & $-4$ & $-4$ & 0 & 0 & 0
			\\ \hline
			$\mathbb{F}_5$ & 0 & 0 & 0 & $\infty$ & $\infty$ & 0 & 0 & 0
			\\ \hline
			$\mathbb{F}_7$ & 0 & 0 & 0 & $-1$ & $\infty$ & 0 & 0 & 0
			\\ \hline
		\end{tabular}
		\label{Tb: FLDR1-help}
	\end{center}
\end{table}
\end{rmk}
\begin{rmk}
	For the values of $H(\Delta)$ for $-200 \le \Delta < 0$,
	see \cite[Table I]{Sch87}.
	We also note that \cite[Proposition 2.2]{Sch87} gives
	a simple formula relating Kronecker's class numbers and
	the class numbers of complex quadratic orders.
	For small $q$ and $n$, we can find a table of
	the values of \Equref{Eq: SchCnt} in \cite{Sch87}.
\end{rmk}

\section{Cubics admitting at most two equivalence classes of
linear determinantal representations}\label{FinLDR}
In this section, we count the number of projective equivalence classes 
of smooth plane cubics over finite fields
admitting at most two equivalence classes of 
linear determinantal representations.

Let $p$ be a prime number, and $m \ge 1$ a positive integer.
Let $\mathbb{F}_q$ be a finite field with $q = p^m$ elements.
Let $\omega \in \mathbb{F}_4$ be an element satisfying
$\omega^2 + \omega + 1=0$.

\begin{thm}\label{Th: 0LDR}\mbox{}
	\begin{enumerate}
		\item\label{0LDR-1} 
		If $q>4$, there are no smooth plane cubics over $\mathbb{F}_q$ 
		which do not admit linear determinantal representations 
		over $\mathbb{F}_q$. 
		\item\label{0LDR-2} 
		If $q \le 4$, there exists only one projective equivalence class of
		smooth plane cubics over $\mathbb{F}_q$ admitting 
		no linear determinantal representations over $\mathbb{F}_q$.
		For explicit representatives of these curves, see Table \ref{Tb: 0LDR}.
	\end{enumerate}
\end{thm}

\begin{prf}
	The assertion follows from Corollary \ref{Th: FinLDR}.
	Here we give another proof of \eqref{0LDR-1}
	which do not use Corollary \ref{Th: FinLDR}.
	Let $C$ be a smooth plane cubic over $\mathbb{F}_q$.
	By the Hasse--Weil bound, we have
	\[
		\#C(\mathbb{F}_q) \ge q + 1 - 2\sqrt{q} = (\sqrt{q} - 1)^2.
	\]
	If $q > 4$, we have $\sqrt{q} > 2$ and
	\[
		\#C(\mathbb{F}_q) > (2-1)^2 = 1.
	\]
	Hence $C$ has at least two $\mathbb{F}_q$-rational points.
	By Proposition \ref{Corr2}, 
	$C$ admits a linear determinantal representation
	over $\mathbb{F}_q$.
\qed\end{prf}

Next, we determine the smooth plane cubics over finite fields
which admit 1 or 2 equivalence classes of
linear determinantal representations.

\begin{thm}\label{Th: 1LDR}\mbox{}
	\begin{enumerate}
		\item\label{1LDR-1} 
		If $q>5$, there are no smooth plane cubics over $\mathbb{F}_q$ 
		admitting a unique equivalence class of 
		linear determinantal representations over $\mathbb{F}_q$. 
		\item\label{1LDR-2} 
		If $q \le 5$, there exists only one projective equivalence class of
		smooth plane cubics over $\mathbb{F}_q$ 
		admitting a unique equivalence
		class of linear determinantal representations over $\mathbb{F}_q$.
		For explicit representatives of these curves, see Table \ref{Tb: 1LDR}.
	\end{enumerate}
\end{thm}

\begin{thm}\label{Th: 2LDR}
	\mbox{}
	\begin{enumerate}
		\item\label{2LDR-1} 
		If $q>7$, there are no smooth plane cubics 
		over $\mathbb{F}_q$ admitting exactly two
		equivalence classes of linear determinantal representations
		over $\mathbb{F}_q$.
		\item\label{2LDR-2} 
		If $q = 2, 3, 5, 7$, there exist 2 projective equivalence classes of
		smooth plane cubics over $\mathbb{F}_q$ admitting exactly two
		equivalence classes of linear determinantal representations
		over $\mathbb{F}_q$.
		\item\label{2LDR-3}
		If $q=4$, there exist 4 projective equivalence classes of smooth plane
		cubics over $\mathbb{F}_q$ admitting exactly two
		equivalence classes of linear determinantal representations
		over $\mathbb{F}_q$.
	\end{enumerate}
	For explicit representatives of the curves in \eqref{2LDR-2} and \eqref{2LDR-3}, 
	see Table \ref{Tb: 2LDR-2} to Table \ref{Tb: 2LDR-7}.
\end{thm}

The proofs of Theorem \ref{Th: 1LDR} and 
Theorem \ref{Th: 2LDR} are omitted 
because they are similar to the proof of Theorem \ref{Th: 0LDR}.

\section{Tables of smooth plane cubics}

Let us show examples of smooth plane cubics
corresponding to cells in Table \ref{Tb: FLDR1}, i.e., smooth plane cubics
over finite fields admitting at most two equivalence classes of
linear determinantal representations.
Moreover, using Theorem \ref{Th: LDR2},
we give a representative of each equivalence class of 
linear determinantal representations of each curve.

Table \ref{Tb: 0LDR} is a summary of smooth plane cubics
over finite fields admitting no linear determinantal representations.

Table \ref{Tb: 1LDR} is a summary of smooth plane cubics over finite fields
admitting a unique equivalence class of linear determinantal representations.

Note that, for these curves in Table \ref{Tb: 1LDR}, 
each linear determinantal representation is
equivalent to a \textit{symmetric} determinantal representation.
For example, in the case of the smooth plane cubic
$X^2Z + XYZ + Y^3 + Y^2Z + YZ^2$ over $\mathbb{F}_2$, 
we transform
\[
\begin{pmatrix}
	0 & 1 & 0 \\
	1 & 0 & 0 \\
	1 & 0 & 1
\end{pmatrix}
\begin{pmatrix}
	0 & Z & Y \\
	Y & 0 & X \\
	X & Y+Z & X + Z
\end{pmatrix} = \begin{pmatrix}
	Y & 0 & X \\
	0 & Z & Y \\
	X & Y & X + Y + Z
\end{pmatrix}.
\]
In fact, symmetric determinantal representations of $C$ are bijective to
$\Pic^0(C)[2] \setminus \{0\}$ (see \cite[Proposition 4.2]{II14}), and
$\Pic^0(C)[2] \cong \Integer / 2\Integer$ for the cubics $C$ 
in Table \ref{Tb: 1LDR}.
By changing the basis $\{e_0, e_1, e_2\}$, 
we have Table \ref{Tb: sym1LDR} of 
symmetric determinantal representations.

Table \ref{Tb: 2LDR-2} to Table \ref{Tb: 2LDR-7} give
summaries of smooth plane cubics over finite fields
admitting exactly two equivalence classes of linear determinantal representations.

\begin{table}[!h]
	\begin{center}
		\vspace{10pt}
		\caption{Smooth plane cubics over finite fields admitting no
		linear determinantal representations.}
		\vspace{-10pt}
		\begin{tabular}{|c|c|c|c|}
			\hline
			$\mathbb{F}_q$ & $F(X, Y, Z)$ & $C(\mathbb{F}_q)$ &
			 $\#\mathrm{LDR}(C)$  \\
			\hline 
			$\mathbb{F}_2$ & $\begin{array}{c}
			X^2Z + XZ^2 + Y^3 + Y^2Z + Z^3
			\end{array}$ & $[1:0:0]$ (flex) & 0 \\
			\hline
			$\mathbb{F}_3$ & $\begin{array}{c}
			X^2Z + Y^3 - YZ^2 + Z^3
			\end{array}$ & $[1:0:0]$ (flex) & 0 \\
			\hline
			$\mathbb{F}_4$ & $\begin{array}{c}
			X^2Z + XZ^2 + Y^3 + \omega Z^3
			\end{array}$ & $[1:0:0]$ (flex) & 0 \\
			\hline
		\end{tabular}
		\label{Tb: 0LDR}
	\end{center}
\end{table}
\begin{table}[!h]
	\begin{center}
		\vspace{10pt}
		\caption{Smooth plane cubics over finite fields admitting a 
		unique equivalence class of linear determinantal representations.}
		\vspace{-10pt}
		\begin{tabular}{|c|c|c|c|c|c|}
			\hline
			$\mathbb{F}_q$ & $F(X, Y, Z)$ & $C(\mathbb{F}_q)$ 
			& $\#\mathrm{LDR}(C)$ & 
			 $\begin{array}{c}
			\mbox{Linear determinantal}\\ 
			\mbox{representations}
			\end{array}$ \\
			\hline 
			$\mathbb{F}_2$ & $\begin{array}{c}
			X^2Z + XYZ + Y^3 \\+ Y^2Z + YZ^2 
			\end{array}$ & 
			\parbox{77pt}{
			$[1: 0: 0]$ (flex), 
			$[0: 0: 1]$}
			&1 &  $\begin{pmatrix}
				0 & Z & Y \\
				Y & 0 & X \\
				X & Y+Z & X + Z
			\end{pmatrix}$ \\
			\hline
			$\mathbb{F}_3$ & $\begin{array}{c}
			X^2Z - Y^3 \\+ Y^2Z + YZ^2 
			\end{array}$ & 
			\parbox{77pt}{
			$[1: 0: 0]$ (flex), 
			$[0: 0: 1]$}
			 & 1 &  $\begin{pmatrix}
				0 & Z & -Y \\
				Y & 0 & -X \\
				X & -Y+Z & Z
			\end{pmatrix}$ \\
			\hline
			$\mathbb{F}_4$ & $\begin{array}{c}
			X^2Z + \omega XYZ + Y^3\\ + Y^2 Z + \omega YZ^2 
			\end{array}$ & 
			\parbox{77pt}{
			$[1:0:0]$ (flex), 
			$[0:0:1]$}
			 &1 & $\begin{pmatrix}
				0 & Z & Y \\
				Y & 0 & X \\
				X & Y+Z & \omega X + \omega Z
			\end{pmatrix}$ \\
			\hline
			$\mathbb{F}_5$ & $\begin{array}{c}
			X^2Z + Y^3 + 2 YZ^2 
			\end{array}$ & 
			\parbox{77pt}{
			$[1:0:0]$ (flex), 
			$[0:0:1]$}
			 &1 & $\begin{pmatrix}
				0 & Z & -Y \\
				Y & 0 & -X \\
				X & Y & 2Z
			\end{pmatrix}$ \\
			\hline
		\end{tabular}
		\label{Tb: 1LDR}
	\end{center}
\end{table}
\begin{table}[!h]
	\begin{center}
		\vspace{10pt}
		\caption{Examples of symmetric determinantal representations for
		smooth plane cubics in Table \ref{Tb: 1LDR}.}
		\vspace{-10pt}
		\begin{tabular}{|c|c|c|c|c|c|}
			\hline
			$\mathbb{F}_q$ & $F(X, Y, Z)$ & 
			 $\begin{array}{c}
			\mbox{Symmetric determinantal}\\ 
			\mbox{representations}
			\end{array}$ \\
			\hline 
			$\mathbb{F}_2$ & $\begin{array}{c}
			X^2Z + XYZ + Y^3 \\+ Y^2Z + YZ^2 
			\end{array}$ & 
			$\begin{pmatrix}
				Y & 0 & X \\
				0 & Z & Y \\
				X & Y & X + Y + Z
			\end{pmatrix}$ \\
			\hline
			$\mathbb{F}_3$ & $\begin{array}{c}
			X^2Z - Y^3 \\+ Y^2Z + YZ^2 
			\end{array}$ & 
			$\begin{pmatrix}
				Y & 0 & -X \\
				0 & -Z & Y \\
				-X & Y & -Y - Z
			\end{pmatrix}$ \\
			\hline
			$\mathbb{F}_4$ & $\begin{array}{c}
			X^2Z + \omega XYZ + Y^3\\ + Y^2 Z + \omega YZ^2 
			\end{array}$ & 
			$\begin{pmatrix}
				Y & 0 & X \\
				0 & Z & Y \\
				X & Y & \omega X + Y + \omega Z
			\end{pmatrix}$ \\
			\hline
			$\mathbb{F}_5$ & $\begin{array}{c}
			X^2Z + Y^3 + 2 YZ^2 
			\end{array}$ & 
			$\begin{pmatrix}
				-Y & 0 & X \\
				0 & -Z & Y \\
				X & Y & 2Z
			\end{pmatrix}$ \\
			\hline
		\end{tabular}
		\label{Tb: sym1LDR}
	\end{center}
\end{table}
\begin{table}[!h]
	\begin{center}
		\vspace{10pt}
		\caption{Smooth plane cubics over 
		$\mathbb{F}_2$ admitting exactly 
		two equivalence classes of linear determinantal representations.}
		\vspace{-10pt}
		\begin{tabular}{|c|c|c|c|c|c|}
			\hline
			$\mathbb{F}_q$ & $F(X, Y, Z)$ & $C(\mathbb{F}_q)$ 
			& $\#\mathrm{LDR}(C)$ & 
			 $\begin{array}{c}
			\mbox{Linear determinantal}\\ 
			\mbox{representations}
			\end{array}$ \\
			\hline 
		$\mathbb{F}_2$ & $\begin{array}{c}
		X^2Z + XY^2 + YZ^2 
		\end{array}$ & 
			\parbox{77pt}{
			$[1:0:0]$, \\$[0:1:0]$, \\
			$[0:0:1]$}
			 & 2 & 
		\(
			\begin{array}{l}
			\begin{pmatrix}
				0 & Z & Y \\
				Z & Y & X \\
				X & 0 & Y
			\end{pmatrix},\\
			\begin{pmatrix}
				0 & Z & Y \\
				Y & 0 & X \\
				X & X & Z
			\end{pmatrix}
			\end{array}
		\) \\
		\cline{2-5}
		 & $\begin{array}{c} 
		 	X^2Z + XZ^2 + Y^3
		 \end{array}$ & \parbox{77pt}{
			$[1:0:0]$ (flex), \\
			$[1:0:1]$ (flex), \\
			$[0:0:1]$ (flex)}
			 & 2 & 
		\(
			\begin{array}{l}
			\begin{pmatrix}
				0 & Z & Y \\
				Y & 0 & X \\
				X + Z & Y & 0
			\end{pmatrix},\\
			\begin{pmatrix}
				0 & Z & Y \\
				Y & 0 & X + Z \\
				X & Y & 0
			\end{pmatrix}
			\end{array}
		\) \\ \hline
		\end{tabular}
		\label{Tb: 2LDR-2}
	\end{center}
\end{table}
\begin{table}[!h]
	\begin{center}
		\vspace{10pt}
		\caption{Smooth plane cubics over $\mathbb{F}_3$ 
		admitting exactly 
		two equivalence classes of linear determinantal representations.}
		\vspace{-10pt}
		\begin{tabular}{|c|c|c|c|c|c|}
			\hline
			$\mathbb{F}_q$ & $F(X, Y, Z)$ & $C(\mathbb{F}_q)$ 
			& $\#\mathrm{LDR}(C)$ & 
			 $\begin{array}{c}
			\mbox{Linear determinantal}\\ 
			\mbox{representations}
			\end{array}$ \\
			\hline 
		%%% 3
		$\mathbb{F}_3$ & $\begin{array}{c}
		X^2Z + XY^2 + YZ^2 + 2XYZ 
		\end{array}$ & \parbox{77pt}{
			$[1:0:0]$, \\$[0:1:0]$, \\
			$[0:0:1]$}
			 & 2 & \(
			\begin{array}{l}
			\begin{pmatrix}
				0 & Z & -Y \\
				Z & Y & X - Y \\
				X & 0 & -Y
			\end{pmatrix}, \\
			\begin{pmatrix}
				0 & Z & -Y \\
				Y & 0 & -X \\
				X & X & -X + Z
			\end{pmatrix}
			\end{array}
		\) \\
		\cline{2-5}
		 & $\begin{array}{c}
		 	X^2Z - XZ^2 - XYZ - Y^3 
		\end{array}$ & \parbox{77pt}{
			$[1:0:0]$ (flex), \\
			$[1:0:1]$ (flex), \\
			$[0:0:1]$ (flex)} & 2 & 
		\(
			\begin{array}{l}
			\begin{pmatrix}
				0 & Z & -Y \\
				Y & 0 & -X \\
				X - Z & -Y & -X
			\end{pmatrix},\\
			\begin{pmatrix}
				0 & Z & -Y \\
				Y & 0 & -X + Z \\
				X & -Y & -X
			\end{pmatrix}
			\end{array}
		\) \\
		\hline
		\end{tabular}
		\label{Tb: 2LDR-3}
	\end{center}
\end{table}
\begin{table}[!h]
	\begin{center}
		\vspace{10pt}
		\caption{Smooth plane cubics over $\mathbb{F}_4$ 
		admitting exactly 
		two equivalence classes of linear determinantal representations.}
		\vspace{-10pt}
		\begin{tabular}{|c|c|c|c|c|c|}
			\hline
			$\mathbb{F}_q$ & $F(X, Y, Z)$ & $C(\mathbb{F}_q)$ 
			& $\#\mathrm{LDR}(C)$ & 
			 $\begin{array}{c}
			\mbox{Linear determinantal}\\ 
			\mbox{representations}
			\end{array}$ \\
			\hline
		%%% 4 %%%
		$\mathbb{F}_4$ & $\begin{array}{c}
		X^2Z + XY^2 + \omega YZ^2 
		\end{array}$ & \parbox{77pt}{
			$[1:0:0]$, \\$[0:1:0]$, \\
			$[0:0:1]$}
			 & 2 & \(
			\begin{array}{l}
			\begin{pmatrix}
				0 & Z & Y \\
				Z & Y & X \\
				X & 0 & \omega Y
			\end{pmatrix}, \\
			\begin{pmatrix}
				0 & Z & Y \\
				Y & 0 & X \\
				X & X & \omega Z
			\end{pmatrix}
			\end{array}
		\) \\
		\cline{2-5}
		 & $\begin{array}{c}
		 	X^2Z + XY^2 + (\omega + 1)YZ^2 
		\end{array}$ & \parbox{77pt}{
			$[1:0:0]$, \\$[0:1:0]$, \\
			$[0:0:1]$}
			 & 2 & 
		\(
			\begin{array}{l}
			\begin{pmatrix}
				0 & Z & Y \\
				Z & Y & X \\
				X & 0 & (\omega + 1) Y
			\end{pmatrix}, \\
			\begin{pmatrix}
				0 & Z & Y \\
				Y & 0 & X \\
				X & X & (\omega + 1) Z
			\end{pmatrix}
			\end{array}
		\) \\
		\cline{2-5}
		 & $\begin{array}{c}
		 	X^2Z + XZ^2 + \omega Y^3 
		\end{array}$ & \parbox{77pt}{
			$[1:0:0]$ (flex), \\
			$[1:0:1]$ (flex), \\
			$[0:0:1]$ (flex)} & 2 & 
		\(
			\begin{array}{l}
			\begin{pmatrix}
				0 & Z & Y \\
				Y & 0 & X \\
				X + Z & \omega Y & 0
			\end{pmatrix},\\
			\begin{pmatrix}
				0 & Z & Y \\
				Y & 0 & X + Z \\
				X & \omega Y & 0
			\end{pmatrix}
			\end{array}
		\) \\
		\cline{2-5}
		 & $\begin{array}{c}
		 	X^2Z + XZ^2 +(\omega+1) Y^3 
		\end{array}$ & \parbox{77pt}{
			$[1:0:0]$ (flex), \\
			$[1:0:1]$ (flex), \\
			$[0:0:1]$ (flex)} & 2 & 
		\(
			\begin{array}{l}
			\begin{pmatrix}
				0 & Z & Y \\
				Y & 0 & X \\
				X + Z & (\omega + 1) Y & 0
			\end{pmatrix},\\
			\begin{pmatrix}
				0 & Z & Y \\
				Y & 0 & X + Z \\
				X & (\omega + 1) Y & 0
			\end{pmatrix}
			\end{array}
		\) \\
		\hline
		\end{tabular}
		\label{Tb: 2LDR-4}
	\end{center}
\end{table}
\begin{table}[!h]
	\begin{center}
		\vspace{10pt}
		\caption{Smooth plane cubics over $\mathbb{F}_5$ 
		admitting exactly 
		two equivalence classes of linear determinantal representations.}
		\vspace{-10pt}
		\begin{tabular}{|c|c|c|c|c|c|}
			\hline
			$\mathbb{F}_q$ & $F(X, Y, Z)$ & $C(\mathbb{F}_q)$ 
			& $\#\mathrm{LDR}(C)$ & 
			 $\begin{array}{c}
			\mbox{Linear determinantal}\\ 
			\mbox{representations}
			\end{array}$ \\
			\hline 
		%%% 5 %%%
		$\mathbb{F}_5$ & $\begin{array}{c}
		X^2Z + XY^2 + YZ^2 -2 XYZ 
		\end{array}$ & \parbox{77pt}{
			$[1:0:0]$, \\$[0:1:0]$, \\
			$[0:0:1]$}
			 & 2 & \(
			\begin{array}{l}
			\begin{pmatrix}
				0 & Z & -Y \\
				Z & Y & X -2Y \\
				X & 0 & -Y
			\end{pmatrix}, \\
			\begin{pmatrix}
				0 & Z & -Y \\
				Y & 0 & -X \\
				X & X & -2X + Z
			\end{pmatrix}
			\end{array}
		\) \\
		\cline{2-5}
		 & $\begin{array}{c}
		 	X^2Z - XZ^2 -2 XYZ - Y^3 
		\end{array}$ & \parbox{77pt}{
			$[1:0:0]$ (flex), \\
			$[1:0:1]$ (flex), \\
			$[0:0:1]$ (flex)}
			 & 2 & 
		\(
			\begin{array}{l}
			\begin{pmatrix}
				0 & Z & -Y \\
				Y & 0 & -X \\
				X - Z & -Y & -2X
			\end{pmatrix},\\
			\begin{pmatrix}
				0 & Z & -Y \\
				Y & 0 & -X + Z \\
				X & -Y & -2X
			\end{pmatrix}
			\end{array}
		\) \\
		\hline
		\end{tabular}
		\label{Tb: 2LDR-5}
	\end{center}
\end{table}

\begin{table}[!h]
	\begin{center}
		\vspace{10pt}
		\caption{Smooth plane cubics over $\mathbb{F}_7$ 
		admitting exactly 
		two equivalence classes of linear determinantal representations.}
		\vspace{-10pt}
		\begin{tabular}{|c|c|c|c|c|c|}
			\hline
			$\mathbb{F}_q$ & $F(X, Y, Z)$ & $C(\mathbb{F}_q)$ 
			& $\#\mathrm{LDR}(C)$ & 
			 $\begin{array}{c}
			\mbox{Linear determinantal}\\ 
			\mbox{representations}
			\end{array}$ \\
			\hline 
		%%% 7 %%%
		$\mathbb{F}_7$ & $\begin{array}{c}
		X^2Z + XY^2 + 3YZ^2 
		\end{array}$ & \parbox{77pt}{
			$[1:0:0]$, \\$[0:1:0]$, \\
			$[0:0:1]$}
			 & 2 & \(
			\begin{array}{l}
			\begin{pmatrix}
				0 & Z & -Y \\
				Z & Y & X \\
				X & 0 & -3Y
			\end{pmatrix}, \\
			\begin{pmatrix}
				0 & Z & -Y \\
				Y & 0 & -X \\
				X & X & 3Z
			\end{pmatrix}
			\end{array}
		\) \\
		\cline{2-5}
		 & $\begin{array}{c}
		 	X^2Z - XZ^2 + 3Y^3 
		\end{array}$ & \parbox{77pt}{
			$[1:0:0]$ (flex), 
			$[1:0:1]$ (flex), \\
			$[0:0:1]$ (flex)}
			 & 2 & 
		\(
			\begin{array}{l}
			\begin{pmatrix}
				0 & Z & -Y \\
				Y & 0 & -X \\
				X - Z & 3Y & 0
			\end{pmatrix},\\
			\begin{pmatrix}
				0 & Z & -Y \\
				Y & 0 & -X + Z \\
				X & 3Y & 0
			\end{pmatrix}
			\end{array}
		\) \\
		\hline\end{tabular}
		\label{Tb: 2LDR-7}
	\end{center}
\end{table}
\clearpage
\subsubsection*{Acknowledgements}
The author would like to thank sincerely Professor Tetsushi Ito 
for various and inspiring comments.

\end{document}